\documentclass[reqno,10pt]{amsart}
\usepackage{latexsym,amssymb,amsthm,amsmath}

\theoremstyle{plain}

\theoremstyle{remark}
\newtheorem{remark}{\bf Remark}

\numberwithin{equation}{section}

\def\<{\left < }
\def\>{\right >}
\def\({\left ( }
\def\){\right )}
\def\e{\eqref}

\def\e{\eqref}

\begin{document}

\title[Open problems and conjectures: recent
development]{Some open problems and conjectures on submanifolds of finite type: recent
development}

\author[B.-Y. Chen]{Bang-yen Chen}
\vskip.15in
\begin{abstract} Submanifolds of finite type were introduced by the author during the late 1970s. The first results on this subject were collected in author's books
\cite{C4,C7}. In 1991, a list of twelve open problems and three conjectures on finite type submanifolds was published in  \cite{C18}.  A detailed survey of the  results,  up to 1996, on this subject was given by the author in \cite{C26}.  Recently, the study of finite type submanifolds, in particular, of biharmonic submanifolds,  have received a growing attention with many progresses since the beginning of this century.  In this article, we provide a detailed account of recent
development on the problems and conjectures listed in \cite{C18}.
 \end{abstract} 

\address{Department of Mathematics \\
	Michigan State University \\East Lansing, Michigan 48824--1027, U.S.A.}
\email{bychen@math.msu.edu}



\maketitle

\vskip.2in
\section{Introduction}

Algebraic geometry studies algebraic varieties which are defined locally as the common zero sets of polynomials.  In algebraic geometry, one can define the degree of an algebraic variety by  its algebraic structure. The concept of degree plays a fundamental role in algebraic geometry. 
On the other hand,  according to  Nash's imbedding theorem, every  Riemannian manifold can be realized as a Riemannian submanifold in some Euclidean space with sufficiently high codimension. However, one lacks the notion of the degree for Riemannian submanifolds in Euclidean spaces. 

Inspired by the above simple observation, the author introduced in the late 1970s  the notions of ``order'' and ``type''  for  submanifolds of  Euclidean spaces and used them to introduce the notion of finite type submanifolds. Just like  minimal submanifolds, submanifolds of finite type can be characterized by a spectral variational principle; namely, as critical points of directional deformations \cite{CDVV2}.

The family of submanifolds of finite type is huge, which contains many important families of submanifolds; including  minimal submanifolds of Euclidean space,  minimal submanifolds of hyperspheres, parallel submanifolds as well as all equivariantly immersed compact homogeneous submanifolds. 

On one hand, the notion of finite type submanifolds  provides a very natural way to apply spectral geometry to study submanifolds.  On the other hand,  one can  also apply the theory of finite type submanifolds to investigate the spectral geometry of submanifolds. For instance,  the author was able to obtain sharp estimates of the total mean curvature for compact submanifolds of Euclidean space via his theory of finite type submanifolds. 

The first results on submanifolds of finite type were collected in \cite{C4,C7}. A list of twelve open problems and three conjectures on submanifolds of finite type was published in 1991  \cite{C18}. Furthermore, a detailed report of the progress  on this theory,  up to 1996,  was presented in \cite{C26}.

Recently, the study of finite type submanifolds, in particular, of biharmonic submanifolds,  have received a growing attention with many progresses since the beginning of this century.  In this article, we provide a detailed account of recent
development on the problems and conjectures listed in \cite{C18}.

\section{Preliminaries}

\subsection{Finite type submanifolds} We recall some basic definitions, results and formulas (for more details see for instance \cite{C4,C30}).

 Let $x:M\rightarrow \mathbb E^m$ be an isometric  immersion of a (connected) Riemannian manifold $M$ into the Euclidean $m$-space $\mathbb E^m$.  Denote by $\Delta$ the Laplace operator of $M$. The immersion $x$ is said to be {\it of finite type\/}  if the position vector field of $M$ in $\mathbb E^m$, also denoted by $x$, can be expressed as a finite sum of $\mathbb E^m$-valued eigenfunctions of the Laplace operator, i.e., if $x$ can be expressed as
\begin{align}\label{2.1}x=c+x_{1}+x_{2}+\ldots+x_{k}\end{align}
where $c$ is a constant vector in $\mathbb E^m$ and $x_{1},\ldots,x_{k}$ are non-constant $\mathbb E^m$-valued maps satisfying \begin{align}\label{2.2}\Delta x_{i}=\lambda_{i}x_{i}, \; i=1,\ldots,k.\end{align} 

The decomposition \e{2.1} is called the {\it spectral decomposition\/} or the {\it spectral resolution} of the immersion $x$.
In particular, if all of the eigenvalues $\lambda_{1},\ldots,\lambda_{k}$ associated with the spectral decomposition are mutually different, then the immersion $x$ (or the submanifold $M$) is said to be of {\it k-type\/}.  In particular, if one of  $\lambda_{1},\ldots,\lambda_{k}$ is zero, then the immersion is said to be of {\it null k-type.}  Clearly,  every submanifold of null $k$-type is non-compact.
A submanifold is said to be {\it of infinite type\/} if it is not of finite type. In terms of finite type submanifolds, a  result of  \cite{Ta1} states that a submanifold of  $\mathbb E^m$ is of 1-type if and only if it is either a minimal
submanifold of $\mathbb E^m$ or a minimal submanifold of a hypersphere of $\mathbb E^m$.

For a spherical isometric immersion $x:M\to S^{m-1}_c\subset \mathbb E^m$,  the immersion is called {\it mass-symmetric\/} in $S_c^{m-1}$ if  the center of gravity of $M$ in $\mathbb E^m$ coincides with the center $c$ of  the hypersphere $S_c^{m-1}$ in $\mathbb E^m$.

\subsection{Minimal polynomials}
For a finite type submanifold $M$ satisfying \e{2.1} and \e{2.2}, the polynomial
$P$ defined by
\begin{align}\label{2.3}P(t)=\prod_{i=1}^{k} (t-\lambda_{i}),\end{align}  
satisfies $P(\Delta)(x-c)=0.$  This polynomial $P$ is called the {\it minimal polynomial\/} of $M$. For an $n$-dimensional submanifold $M$ of a Euclidean space, the mean curvature vector $H$ satisfies Beltrami's formula:
\begin{align}\label{2.4}\Delta x = -nH.\end{align}

It follows from \e{2.4}  that the minimal polynomial $Q$  also satisfies  $Q(\Delta)H=0$.  Conversely,  if $M$ is compact and if there exists a constant vector $c$ and a nontrivial polynomial $Q$ such that $Q(\Delta)(x-c)=0$ (or $Q(\Delta)H=0$), then $M$ is always of finite type.  This characterization of finite type submanifolds via the minimal polynomial plays an important role in the theory of finite type submanifolds. 

When $M$ is non-compact,  the existence of a nontrivial polynomial $Q$  satisfying $Q(\Delta)H=0$ does not guarantee $M$ to be finite type. On the other hand, if either $M$ is one-dimensional or $Q$ is a polynomial  of degree $k$ with exactly $k$ distinct  roots,  then the existence of the polynomial $Q$ satisfying $Q(\Delta)(x-c)=0$ for some constant vector $c$ does guarantee that $M$ is of finite type  \cite{CP}.

\subsection{A basic formula for $\Delta H$}

The following basic formula of $\Delta H$ derived in \cite{C2,C4,C7} plays important role in the study of submanifolds of low type as well as in the study of biharmonic submanifolds: 
\begin{align}\label{2.5} \Delta H=\Delta^{D}H + \sum_{i=1}^{n} h(e_{i},A_{H}e_{i}) + 2\,{\rm trace}\,(A_{DH})+ {n\over 2} {\rm grad} \<H,H\> ,\end{align} 
where $\Delta^{D}$ is the Laplace operator associated with the normal connection $D$, $h$ is the second fundamental form, and $\{ e_{1},\ldots,e_{n}\}$ is a local orthonormal frame of $M$. In particular, if $M$ is a hypersurface of a Euclidean space $\mathbb E^{n+1}$, then formula \e{2.5} reduces to
\begin{align}\label{2.6} \Delta H=(\Delta \alpha +\alpha ||h||^{2})\xi + 2\, {\rm trace}\,(A_{DH})+{n\over 2}\, {\rm grad}\<H,H\>,\end{align}
where $\alpha$ is the mean curvature and $\xi$ a unit normal vector of $M$ in $\mathbb E^{n+1}$.

Similar formulas hold as well if the ambient spaces is pseudo-Euclidean.

\subsection{$\delta$-invariants and ideal immersions}  Let $M$ be a  Riemannian $n$-manifold. Denote by $K(\pi)$ the sectional curvature of  a plane section $\pi\subset T_pM$, $p\in M$. For any orthonormal basis $e_1,\ldots,e_n$ of $T_pM$, the scalar curvature $\tau$ at $p$ is
 $$\tau(p)=\sum_{i<j} K(e_i\wedge e_j).$$

Let $L$ be a $r$-subspace of $T_pM$ with $r\geq 2$  and let $\{e_1,\ldots,e_r\}$ be an orthonormal basis of $L$. The scalar curvature $\tau(L)$ of  $L$ is defined by
\begin{align}\label{2.8}\tau(L)=\sum_{\alpha<\beta} K(e_\alpha\wedge e_\beta),\quad 1\leq \alpha,\beta\leq r.\end{align}

For given integers $n\geq 3$,  $k\geq 1$, we denote by $\mathcal S(n,k)$ the finite set  consisting of $k$-tuples $(n_1,\ldots,n_k)$ of integers  satisfying  
$2\leq n_1,\cdots,
n_k<n$  and $\sum_{{j=1}}^{k} n_i\leq n.$ 

Put ${\mathcal S}(n)=\cup_{k\geq 1}\mathcal S(n,k)$. For each $k$-tuple $(n_1,\ldots,n_k)\in \mathcal S(n)$, the author introduced in 1990s the Riemannian invariant $\delta{(n_1,\ldots,n_k)}$ by
\begin{align} \delta(n_1,\ldots,n_k)(p)=\tau(p)- \inf\{\tau(L_1)+\cdots+\tau(L_k)\},\;\; p\in M,\end{align} where $L_1,\ldots,L_k$ run over all $k$ mutually orthogonal subspaces of $T_pM$ such that $\dim L_j=n_j,\, j=1,\ldots,k$ (cf. \cite{C29} for details).

For an $n$-dimensional submanifold of $\mathbb E^m$ and for a $k$-tuple  $(n_1,\ldots,n_k)\in\mathcal S(n)$,  the author proved the following general sharp inequality  \cite{C26,C30}:
\begin{align}\label{2.9} \delta{(n_1,\ldots,n_k)} \leq \text{$ {{n^2(n+k-1-\sum n_j)}\over{2(n+k-\sum n_j)}}$} |H|^2 ,\end{align}
where $|H|^2=\<H,H\>$ denotes the squared mean curvature of $M$. 

A submanifold $M$ of $\mathbb E^m$ is called {\it $\delta(n_1,\ldots,n_k)$-ideal} if it satisfies the equality case of \eqref{2.9} identically. 
  Roughly speaking, ideal submanifolds are submanifolds which receive the least possible tension from its ambient space. For the most recent survey on $\delta$-invariants and ideal immersions, see \cite{C30,C31} for details.

\subsection{Proper and $\epsilon$-superbiharmonic submanifolds}
An immersed submanifold $M$ of a Riemannian manifold $\tilde M$ is said to be {\it properly
immersed} if the immersion is a proper map, i.e.,  the preimage of each compact set in
$\tilde M$ is compact in $M$. 

A hypersurface of a Euclidean space is called {\it weakly convex} if it has non-negative principle curvatures. Also, a hypersphere of an $(n+1)$-sphere is called  isoparametric if it has constant principal curvatures.

The {\it total mean curvature} of a submanifold $M$ in a Riemannian manifold is defined to be $\int_M | H|^2 dv$.

 Let $M$ be a submanifold of a Riemannian manifold  with inner product $\left<\;\,,\;\right>$. Then $M$ is called {\it $\epsilon$-superbiharmonic} if it satisfies
\begin{align}\left<\right.\!\Delta  H, H\! \left.\right>\geq (\epsilon-1)|\nabla H|^2,\end{align}
where $\epsilon\in [0,1]$ is a constant.

For a complete Riemannian manifold $(N,h)$ and $\alpha\geq 0$, if the sectional curvature $K^N$ of $N$ satisfies
\begin{align}K^N\geq -L(1+{\rm dist}_N(\,\cdot\,,q_0)^2)^{{\alpha}/{2}}\end{align}
for some $L>0$ and $q_0\in N$, 
then we  call that $K^N$ has a polynomial growth bound of order $\alpha$ from below.

\vskip.4in 
\section{Finite Type Hypersurfaces of Euclidean Space}

The class of finite type submanifolds in Euclidean spaces is huge. It includes all minimal submanifolds of Euclidean spaces, all minimal submanifolds of  hyperspheres as well as all compact homogeneous submanifolds equivariantly immersed in some Euclidean space. In contrast, very few examples of finite type hypersurfaces in Euclidean spaces are known. So far,  minimal surfaces,  circular cylinders and the spheres are the only known surfaces of finite type in $\mathbb E^3$. 

In \cite{C18}, the author asked the following problem. 
\vskip.1in

 {\bf Problem 1.} {\it Classify all finite type hypersurfaces in $\mathbb E^{n+1}$. In particular, classify all finite type surfaces in $\mathbb E^3$.}
 \vskip.1in
 
This problem does have a complete solution when $n=1$. In fact, it was proved in \cite{C3,C4} that circles are the only finite type closed planar curves.  Also, it is known that lines are the only non-closed planar curves  \cite{C14}. In fact, lines are the only null finite type planar curves  \cite{CDVV}.  
\vskip.1in

Next, we recall some classical results concerning Problem 1. The first result on the classification of finite type surfaces in $\mathbb E^3$ was obtained  in \cite{C11}  which states that  circular cylinders are the only tubes  of finite type.  It was proved  in \cite{CDVV} that a ruled surface in $\mathbb E^3$ is of finite type if and only if it is a plane, a circular cylinder or a helicoid. Furthermore, it was shown in \cite{CD2} that  spheres and  circular cylinders are the only quadrics of finite type in $\mathbb E^3$. In \cite{G2}, it was proved that a cone in $\mathbb E^3$ is of finite type if and only if it is a plane.  
Also, it was shown in \cite{CLu} that every compact 2-type hypersurface of $\mathbb E^{n+1}$ has non-constant mean curvature.

For compact finite type surfaces in $\mathbb E^3$, the author made the following conjecture in
\cite{C7,C18}.
\vskip.1in

 {\bf Conjecture 1.} {\it The only compact surfaces of finite type in $\mathbb E^3$ are the ordinary spheres.}
\vskip.1in

Besides the classical results given above,  there are several additional results  obtained in \cite{AGM,B,C11,DDV1,DPV,HV4,V2} which support this conjecture. However, after more than two decades Conjecture 1 remains open.

In addition to Conjecture 1, the author would like to make the following two additional conjectures which are closely related to Conjecture 1.

\vskip.1in

 {\bf Conjecture 1.A.} {\it The only surfaces of finite type in $\mathbb E^3$ are minimal surfaces, and open portions of spheres and circular cylinders.} 

\vskip.1in

 {\bf Conjecture 1.B.} {\it The only compact  hypersurfaces of finite type in Euclidean space are ordinary hyperspheres.}

\section{Spherical hypersurfaces of finite type}

\subsection{Finite type spherical hypersurfaces}
 In contrast to Euclidean hypersurfaces, there do exist many 1-type and 2-type spherical hypersurfaces. The author proved in \cite{C8} that every compact hypersurface of a hypersphere $S^{n+1}\subset \mathbb E^{n+2}$, not a small hypersphere, is mass-symmetric and of 2-type if and only if it has non-zero constant mean curvature and constant scalar curvature. Consequently,  every isoparametric hypersurface of a hypersphere is either of 1-type or mass-symmetric and 2-type. Since there exist non-minimal isoparametric hypersurfaces in hyperspheres, there do exist  2-type hypersurfaces of hyperspheres.  

The following problem was proposed by the author in \cite{C18}. 
  \vskip.1in

{\bf Problem 2.} {\it Study and classify 2-type hypersurfaces in a hypersphere of $\mathbb E^{n+2}$. In particular, classify 2-type hypersurfaces of a hypersphere $S^4$ in $\mathbb E^5$.}

\vskip.1in

 It is known that a compact surface $M$ in $S^3$ is of 2-type if and only if it is the product of two plane circles of different radii, i.e., $M= S^1(a)\times S^1(b)$ with $a\ne b$ and $ a^2+b^2=1$  \cite{BG,C3}. The same result also holds without compactness \cite{HV1}. 

For 2-type hypersurfaces in $S^4$, we have the following classification theorem from \cite{C24}.
A compact hypersurface of $S^{4}(1)$ is  of $2$-type if and only if it is congruent to one of the following two hypersurfaces:

(a) a standard imbedding $S^1\times S^2\subset S^4(1)\subset \mathbb E^5$  such that the radii $r_1$ of $S^1$  and $r_2$ of $S^2$ satisfying $r_1^2+r_2^2=1$ and $(r_1,r_2)\ne \big(\frac{1}{\sqrt{3}},\frac{\sqrt{2}}{\sqrt{3}}\,\big)$;  

(b) a tube $T^r(V^2)$ with radius $r\ne \frac{\pi}{2}$ over the Veronese surface $V^2$ in $S^4(1)$.

  It was proved by Hasanis and Vlachos in \cite{HV2}  that every 2-type hypersurface of a hypersphere $S^{n+1}$ has nonzero constant mean curvature and constant scalar curvature. 

The following  problems was also proposed by the author in \cite{C18}.
\vskip.1in

{\bf Problem 3.} {\it Classify finite type hypersurfaces of a hypersphere in $\mathbb E^{n+2}$.}
\vskip.1in

  Very little were known on finite type spherical hypersurfaces with type number $\geq 3$.  The only known general result in this respect is that  every 3-type spherical hypersurface  has non-constant mean curvature \cite{C15,CL}.  For finite type spherical surfaces in a 3-sphere $S^3$ with arbitrary type number, the author and Dillen proved in \cite{CD1} that standard 2-spheres and  products of two plane circles are the only compact finite type  surfaces with constant Gauss curvature in $S^3$. 

\vskip.1in

The classification of $k$-type spherical hypersurfaces with $k\geq 3$ remains a very challenge problem. 
From all available information, it seems to the author that there exist no surfaces of $k$-type in $S^3$ for any finite $k$ greater than 2. 
 
\vskip.1in

At an international conferences held at Berlin in 1990, the author proposed the following conjecture.\vskip.1in

 {\bf Conjecture 2.} {\it Minimal surfaces, standard 2-spheres and products of two plane circles are the
only finite type surfaces in $S^3\subset \mathbb E^4$.}
\vskip.1in

This conjecture was also proposed in \cite{C18}.
There are several results obtained in \cite{C15,CBG,CL,FL,Na} which support this conjecture.
 On the other hand, this conjecture remains open in general.

\subsection{Dupin hypersurfaces}

A hypersurface $M$ of $S^{n+1}$ is called a {\it Dupin hypersurface} if the multiplicities of the principal curvatures are constant and each principal curvature are constant along its principal directions.  Since 2-type  spherical  hypersurfaces are mass-symmetric \cite{HV2}, a result of \cite{C8} implies that if a compact Dupin hypersurface $M$ of $S^{n+1}$ is not of 1-type and if it has at most three distinct principal curvatures, then it is parametric if and only if it is of 2-type.

 For Dupin hypersurfaces, the author gave the following open problem in \cite{C18}.
\vskip.1in

{\bf Problem 4.} {\it When is a Dupin hypersurface of  a hypersphere to be of finite type? When is a finite type Dupin hypersurface of a hypersphere to be isoparametric?}
\vskip.1in

As far as the author know, no further results were known for this problem.

\section{Spherical 2-type submanifolds of higher codimension}

It is well-known that there exist ample examples of 1-type surfaces lying fully in odd-dimensional spheres as well as in even-dimensional spheres.  Also, it is known that there exist abundant examples of mass-symmetric 2-type surfaces lying fully in  odd-dimensional hyperspheres (see,  e.g., \cite{BC1,C4, G1, K, Mi}). In contrast,  it was proved by the author and Barros in \cite{BC1} that there do not exist mass-symmetric 2-type surfaces lying fully in $S^4$. So far there are no known examples of non-mass-symmetric 2-type surfaces in $S^4$. In this respect,  the author   asked in \cite{C18} the following 
\vskip.1in

{\bf Problem 5.} {\it Do there exist non-mass-symmetric 2-type surfaces in $S^4$?} 
\vskip.1in

Since there exist no known examples of mass-symmetric 2-type surfaces lying fully in a  hypersphere of a Euclidean space for any even codimension, the author proposed the following problem  in \cite{C18} which is more general than Problem 5.
\vskip.1in

{\bf Problem 6.} {\it Do there exist  2-type surfaces lying fully in an even-dimensional
hypersphere of a Euclidean space\,? In particular, do there exist mass-symmetric 2-type surfaces  lying fully in an even-dimensional hypersphere\,?}
\vskip.1in

If the mass-symmetric, 2-type,  spherical surfaces come from one of the following four families

\begin{itemize}
\item stationary surfaces \cite{BC1}, 

\item  topological 2-spheres\cite{K}, 

\item surfaces with constant Gauss curvature \cite{Mi},  

\item flat Chen surfaces \cite{G1},
\end{itemize}
the answers to Problem 6 is known to be negative.

Up to author's knowledge, Problems 5 and 6 remain unsolved.

\section{Linearly Independent  submanifolds}

The notion of linearly independent immersions or submanifolds  were introduced by the author in \cite{C17}. Suppose that $x:M\rightarrow \mathbb E^m$ is a $k$-type immersion  whose  spectral decomposition is given by \e{2.1}.  Denote by $E_{i}$  the subspace of $\mathbb E^m$ spanned by $\{ x_{i}(p),\, p \in M\}, \,\,i\in \{1,\ldots,k\}$.  The immersion $x$ (or the submanifold $M$)  is called  {\it linearly independent\/} if the subspaces $E_{1},\ldots,E_{k}$ are linearly independent. The
immersion $x$  is called {\it orthogonal} if the subspaces $E_{1},\ldots,E_{k}$ are mutually orthogonal. 

Clearly, every orthogonal immersion is linearly independent and every 1-type immersion is orthogonal.  There exist many examples of  orthogonal immersions and abundant
examples of linearly independent immersions which are not orthogonal. In fact,  every $k$-type curve lying fully in $\mathbb E^{2k}$ and every null $k$-type curve lying fully in $\mathbb E^{2k-1}$ are linearly independent curves; but  $W$-curves are the only orthogonal curves in a Euclidean space. 

  For a given linearly independent immersion $x:M\rightarrow \mathbb E^m$ and a given point $p\in M$, one has the notion of the {\it adjoint hyperquadric\/} $Q_p$  defined in \cite{C17}. When 
 $M$ lies in one of the adjoint hyperquadrics $Q_p,$  $p\in M$, then all of the adjoint hyperquadrics $Q_p,\, p\in M,$ coincide. This common hyperquadric  is  called the {\it adjoint hyperquadric of the linearly independent immersion} \cite{C17}.

It was shown in \cite{C17} that if  $x: M\rightarrow \mathbb E^m$ is  a linearly independent immersion  of a compact manifold, then the submanifold lies in its adjoint hyperquadric  if and only if the submanifold  is spherical with an appropriate center. Moreover, it is known that  a non-minimal linearly independent immersion $x$ is 
orthogonal  if and only if  $M$ is immersed by $x$ as a minimal submanifold of the adjoint hyperquadric  \cite{C17}. Consequently,  every orthogonal immersion of a compact manifold  is  spherical. Moreover, it also known that each  compact homogeneous submanifold, equivariantly immersed in $\mathbb E^m$,  is orthogonal and therefore it is immersed as a minimal submanifold in its adjoint hyperquadric  \cite{C17}.  

It is known that the only linearly independent Euclidean hypersurfaces are hyperspheres, minimal hypersurfaces or spherical hypercylinders (see \cite{CP, DPV,G3,HV4, CDVV2}). Also, by applying  the   classification theorem of 2-type curves  in Euclidean space from \cite{CDV}, one may  conclude that the only linearly independent curves of codimension 2 in a Euclidean space are circles, lines and circular helices. In this respect, the author would like to point out that there do exist abundant examples of
linearly independent curves of codimension 3 in Euclidean space.

In \cite {C18} the next two problems on linearly independent immersions were proposed.
\vskip.1in

{\bf Problem 7.} {\it Study and classify linearly independent  2-type immersions.}
\vskip.1in

{\bf Problem 8.} {\it Study and classify linearly independent  submanifolds of codimension 2.}
\vskip.1in

In  \cite{Ja}, Jang studied  linearly independent immersions {\it with codimension} $\leq 3$. He derived some necessary and sufficient conditions for linearly independent immersions with codimension $\leq 3$ to be orthogonal. His results provide some partial generalizations of author's results in \cite{C17} for codimension $\leq 3$.

The class of linearly independent immersions lies in a much larger class of immersions; namely, the class of  {\it immersions of restricted type} introduced in \cite{CDVV3}. 

A submanifold of a Euclidean  space is said to be of restricted type if its shape operator with respect to the mean curvature vector is the restriction of a fixed linear transformation of the ambient space to the tangent space of the submanifold at every point of the submanifold.  There are very few known results on submanifolds of restricted type. The only known classification results for submanifolds of restricted in Euclidean spaces are the classifications of planar curves and Euclidean hypersurface of restricted type \cite{CDVV3} (see also \cite{V1,V2}).

\vskip.2in

\section{Null 2-type submanifolds}

It was proved in  \cite{CL}  that  every  2-type submanifold in a Euclidean space with parallel mean curvature vector is either spherical  or of null 2-type.  Related with this  the author 
also asked in \cite{C18} the following.
\vskip.1in

{\bf Problem 9.} {\it Is every $n$-dimensional  non-null 2-type submanifold with constant mean curvature be in $\mathbb E^{n+2}$    spherical?}
\vskip.1in

It follows from the definition of null 2-type submanifolds and formula \e{2.3} that the mean curvature vector $H$ of a null  2-type submanifold satisfies \begin{equation}\label{8.1}\Delta H=\lambda H,\end{equation}
where  $\lambda$ is a nonzero real number.   It was proved in \cite{C14} that biharmonic submanifolds, null 2-type submanifolds and 1-type submanifolds are the only Euclidean submanifolds satisfying  \e{8.1}.  

From the classification of finite type planar curves, we know that there do not exist null 2-type planar curves. Furthermore,  it follows from \cite{C14} that the only null 2-type curves in
Euclidean spaces are circular helices  with nonzero torsion in $\mathbb E^3$.  A general result 
 from \cite{C13} states that circular cylinders are the only null 2 type surfaces in Euclidean 3-space. 
Since biharmonic and null 2-type surfaces in $\mathbb E^3$ were classified, the complete classification of surfaces in $\mathbb E^3$  satisfying \e{8.1} was done.   

Next, we present some classical results on null 2-type Euclidean  surfaces with codimension $\geq 2$. It is known in \cite{C14} that  null 2-type surfaces in $\mathbb E^4$ are helical cylinders if they have constant mean curvature.  However, it still unknown whether every null 2-type surface in $\mathbb E^4$ has constant mean curvature.  

 As a generalization of \cite{C13},  null 2-type hypersurfaces with at most two distinct principal curvatures were classified in \cite{FL}   (see \cite{FGL2} for null 2-type conformally flat hypersurfaces of dimension $\not= 3$). In \cite{C14}, the author proved that a surface $M$ in $\mathbb E^4$ is of null 2-type with parallel normalized mean curvature vector if and only if it is an open portion of a circular cylinder in a hyperplane of $\mathbb E^4$. The author also proved that the only null 2-type surfaces in $\mathbb E^4$ with constant mean curvature are open portion of helical cylinders. Hasani and Vlachos proved in \cite{HV6} that non-spherical hypersurfaces in $\mathbb E^4$ with non-vanishing constant mean curvature and constant scalar curvature are the only null 2-type hypersurfaces.

For null 2-type submanifolds, the author proposed in \cite{C18} the following. 
\vskip.1in

{\bf Problem 10.} {\it Study and classify null 2-type submanifolds. In particular, classify all null 2-type
surfaces in 4-dimensional Euclidean space and in 4-dimensional pseudo-Euclidean spaces.}
\vskip.1in

Now, we present some later development concerning Problem 10. It was showed by Li in \cite{L1}  that a surface  in $\mathbb E^m$ with parallel normalized mean curvature vector is of null 2-type if and only if it is an open portion of a circular cylinder. Also, it was proved in \cite{L2} that, for a non-pseudo-umbilical Chen surface $M$ in $\mathbb E^m$,  if $M$ is of null 2-type and with constant mean curvature, then $M$ is flat and must lie  in a  totally geodesic $\mathbb E^6\subset \mathbb E^m$.

In \cite{Du1}, Dursun classified  3-dimensional null 2-type submanifold  of  $\mathbb E^5$ with two distinct principal curvatures in the parallel mean curvature direction and with constant squared norm of the second fundamental form.  
 Also, he proved in \cite{Du3} that if a null 2-type submanifold of $\mathbb E^{n+2}$ with codimension 2  has flat normal connection, constant mean curvature and non-parallel mean curvature vector, then the first normal space must one-dimensional. By using this fact, he derived some some classification results.

Recently, it was proved by the author and Garay in \cite{CG} that  a null 2-type hypersurface of $\mathbb E^{n+1}$ is $\delta(2)$-ideal if and only if it an open portion of a spherical cylinder $S^{n-1}\times \mathbb E$.

However, until now there are still no complete classification of null 2-type submanifolds. In particular, null 2-type surfaces in Euclidean 4-space $\mathbb E^4$ are not completely classified.

\section{Finite Type Submanifolds in Homogeneous Spaces} \vskip.05in

For finite type submanifolds of compact irreducible homogeneous manifolds,  the following problem was proposed in \cite{C18}.
\vskip.1in

{\bf Problem 11.} {\it Let $\tilde M$ be a compact irreducible homogeneous manifold immersed in a Euclidean space $\mathbb E^N$ by its first standard immersion $\phi$ and $M$ a
submanifold of $\tilde M$. When $M$ is of finite type in $\mathbb E^N$ via $\phi$? In particular, when $M$ is of 1- or 2-type in $\mathbb E^N$ via $\phi$?}
\vskip.1in

 If $\tilde  M$ is a projective $m$-space $FP^m$ over a field $F={\mathbb R}, {\mathbb C}$ or  ${\mathbb H}$ with a standard Riemannian metric,  this problem has been investigated in \cite{BC2,BU,C4,C5,MR,UD1}, \cite{Ros1}-\cite{Ros3} and \cite{Di3}-\cite{Di8}, among others.  When $\tilde M$ is  the real Grassmannian $G^{R}(p,q)$ or the space $U(n)/O(n)$, this problem has been investigated in \cite{BN}. 
 
Several recent results in this respect were given by Dimitri\'c  in a series of his papers  \cite{Di5}-\cite{Di8}. In \cite{Di8} he studied  2 and 3-type Hopf hypersurfaces of complex projective space $CP^m$ and of complex hyperbolic space $CH^m$ via some suitable imbeddings into pseudo-Euclidean spaces of Hermitian matrices. He proved in \cite{Di8}  that tubes of certain radii around totally geodesic $CP^k\subset CP^m$, $k\in \{0,\ldots,m-1\}$, and around the complex quadric $Q^{m-1}\subset CP^m$ are 2-type Hopf hypersurfaces in $CP^m$. Conversely, every 2-type Hopf hypersurface in $CP^m$ is locally congruent to such a tube. For $CH^m$, he shown that a Hopf hypersurface is of 2-type if and only if it is locally congruent to a geodesic hypersphere or to a tube  of any radius $r>0$ around a totally geodesic $CH^{m-1}\subset CH^m$. He also obtained  partial classification of 3-type Hopf hypersurfaces in $CP^2$ as well as in $CH^2$.

\section{Biharmonic submanifolds}
The study of biharmonic submanifolds was initiated in the middle of 1980s  via author's study of finite type submanifolds; also independently by Jiang  \cite{J} via his study of Euler-Lagrange's equation of bienergy functional in the sense of Eells and Lemaire \cite{EL1,EL2}. 

Let $x: M \rightarrow \mathbb E^m$ be an isometric immersion. As we already mentioned earlier, the position vector field of $M$ in $\mathbb E^m$ satisfies Beltrami's formula:
\begin{align}\label{10.1}\Delta x=-nH. \end{align}
   Formula \e{10.1} implies that the immersion is minimal if and only if it is harmonic, i.e., $\Delta x = 0.$ An immersion $x : M \rightarrow \mathbb  E^m$ is called {\it biharmonic\/}  if 
   \begin{align}\label{10.2}\Delta^{2}x=0,\hskip.2in or\, equivalently  \hskip.2in \Delta H=0,\end{align}
holds identically. 

Let  $x :M\to \mathbb E^m$ be an isometric immersion.  It follows from \e{2.5} and \e{10.2} that $M$ is a biharmonic submanifold if and only if it satisfies the following fourth order strongly elliptic semi-linear PDE system (see, e.g., \cite{C2,C30}):
 \begin{align} & \Delta^{D} H + \sum_{i=1}^{n}
\sigma(A_{H}e_{i},e_{i})=0,\\ &n\,{\rm grad}\!\left<\right.\! H, H\! \left.\right> + 4\, {\rm trace}\,A_{D H}=0,\end{align}
where $\{e_1,\ldots,e_n\}$ is an orthonormal frame of $M$.

It is obvious that  minimal immersions are trivially biharmonic. Thus, the real problem is if there are other submanifolds besides minimal ones that are biharmonic.

 In \cite{C18}  the author made the following simple geometric question. 
\vskip.1in

{\bf Problem 12.} {\it Other than minimal submanifolds, which submanifolds of $\mathbb E^m$ are biharmonic?}
\vskip.1in

A {\it biharmonic map} is a map $\phi:(M,g)\to (N,h)$ between Riemannian manifolds that is a critical point of the bienergy functional:
\begin{align}\label{10.5} E^2(\phi,D)=\frac{1}{2}\int_D ||\tau_\phi||^2* 1\end{align}
for every compact subset $D$ of $M$, where $\tau_\phi={\rm trace}_g\nabla d\phi$ is the tension field $\phi$.   

The Euler-Lagrange equation of \e{10.5} gives the following biharmonic map
equation  \cite{J}:
\begin{align}\label{BH} \tau^2_\phi:={\rm trace}_g(\nabla^\phi\nabla^\phi-\nabla^\phi_{\nabla^M})\tau_\phi-{\rm trace}_g R^N(d\phi,\tau_\phi)d\phi=0,\end{align}
where $R^N$ denotes the curvature tensor of $(N,h)$. Equation \eqref{BH} implies that $\phi$ is a biharmonic map if and only if its bi-tension field $\tau^2_\phi$ vanishes.  

For an $n$-dimensional submanifold $M$ of $\mathbb E^m$, if we denote by $\iota:M\to \mathbb E^m$ the inclusion map of the submanifold, then the tension field of the inclusion map is given by $\tau_\iota=-\Delta\iota=-n H$ according to Beltrami's formula. Therefore, the submanifold $M$ is biharmonic if and only if 
$$n\Delta  H=- \Delta^2 \iota=-\tau^2_\iota=0,$$
i.e., the inclusion map $\iota$ is a biharmonic map.

\section{Biharmonic conjectures}
\subsection{The original biharmonic conjecture}

The author shown in 1985 that biharmonic surfaces in $\mathbb E^3$ are minimal (independently by Jiang \cite{J}). This result was the starting point of Dimitri\'c's work on his thesis   \cite{Di1} at Michigan State University.   In \cite{Di1}, Dimitri\'c extended author's result (unpublished then) to  biharmonic Euclidean hypersurfaces with at most  two distinct principal curvatures \cite{Di1}.  He also proved that each biharmonic submanifold of finite type  in in a Euclidean space is minimal. Furthermore, he proved that  biharmonic Euclidean submanifolds of finite type and pseudo-umbilical biharmonic Euclidean  submanifolds are always minimal. 

In \cite{C18}, the author pointed out that spherical biharmonic submanifolds of a Euclidean space are minimal. Moreover, it was  proved by Hasani and Vlachos in \cite{HV5} that biharmonic hypersurfaces of ${\mathbb{E}}^4$ are also minimal. 

In \cite{C18}, the author proposed  the following Biharmonic Conjecture.
\vskip.1in

 {\bf Conjecture 3}: \emph{The only biharmonic submanifolds of Euclidean spaces are 
the minimal ones.}
\vskip.1in

\subsection{Generalized Chen's biharmonic conjectures}

 Caddeo,  Montaldo and  Oniciuc proved in \cite{CMO02} that every biharmonic surface in the hyperbolic $3$-space $H^3$ is minimal. They also proved in \cite{CMO01} that biharmonic hypersurfaces of $H^n$ with at most two distinct principal curvatures are minimal.  Based on these facts, they made the following  conjecture  in  \cite{CMO01}.
\vskip.1in

{\bf Generalized Chen's Conjecture}: 
\emph{Every biharmonic submanifold of a Riemannian manifold with non-positive sectional curvature is minimal.}
\vskip.1in

The study of  biharmonic submanifolds is nowadays a very active research subject. Biharmonic submanifolds have received a growing attention with many progresses done since the beginning of this century.

\section{Recent developments  on original biharmonic conjecture}

From the definition of biharmonic submanifolds and Hopf's lemma we see that biharmonic submanifolds in a Euclidean space are always non-compact.

\vskip.06in
The following recent results provides strong supports to the original biharmonic conjecture.

\begin{itemize}
\item Biharmonic properly immersed  \cite{M12a}.

\item Biharmonic  submanifolds which are complete and proper   \cite{AM}.

\item $\delta(2)$-ideal and $\delta(3)$-ideal biharmonic  hypersurfaces  \cite{CM}.

\item Weakly convex biharmonic submanifolds  \cite{Luo1}.

\item Submanifolds whose $L^p,\, p\geq 2$,  integral of the mean curvature vector field satisfies
certain decay condition at infinity  \cite{Luo3}.

\item Biharmonic  submanifolds  satisfying the
decay condition at infinity $$\lim_{\rho\to \infty}\frac{1}{\rho^2}\int_{f^{-1}(B_\rho)}| H |^{2}dv=0,$$
where $f$ is the immersion, $B_\rho$ is a geodesic ball of $N$ with radius $\rho$  \cite{Wh}.

\end{itemize}

The author would like to point out that  Y.-L. Ou showed recently in \cite{Ou09} that the original biharmonic conjecture cannot be generalized to biharmonic conformal submanifolds in Euclidean spaces. 
\vskip.05in

\begin{remark}  Conjecture 3 remains open.\end{remark}

\begin{remark}  Conjecture 3 is false if the ambient Euclidean space were replaced by a pseudo-Euclidean space. The simplest examples are constructed  in \cite{CI91}. \end{remark}

\section{Recent developments  on generalized biharmonic conjecture}

In the last few years, there are many partial answers support the generalized Chen's biharmonic conjecture. The following is a list of recent results which support the generalized biharmonic conjecture.

\begin{itemize}
\item  Biharmonic hypersurfaces in $H^4(-1)$  \cite{BMO10b}.

\item Pseudo-umbilical biharmonic submanifolds of $H^m(-1)$ \cite{CMO01}.

\item Totally umbilical biharmonic hypersurfaces in Einstein spaces \cite{Ou10}.

\item Biharmonic hypersurfaces with finite total mean curvature in a Riemannian manifold of non-positive Ricci curvature  \cite{NU1}.

\item Biharmonic submanifolds with finite total mean curvature in a Riemannian manifold of non-positive sectional curvature  \cite{NU2}.

\item Biharmonic properly immersed submanifolds in a complete Riemannian manifold with non-positive sectional curvature whose sectional curvature has polynomial growth bound of order less than 2 from below  \cite{M12b}.

\item Complete biharmonic submanifolds  with finite bi-energy and energy in a non-positively curved Riemannian manifold \cite{NUG}.

\item Complete oriented biharmonic hypersurfaces $M$  whose mean curvature $H$ satisfying $H\in L^2(M)$  in a Riemannian manifold with non-positive Ricci tensor  \cite{Al}.

\item Compact biharmonic submanifolds in a Riemannian manifold with non-positive sectional curvature  \cite{M13}.

\item Complete biharmonic submanifolds (resp., hypersurfaces) in a Riemannian manifold
whose sectional curvature (resp., Ricci curvature) is non-positive with at most polynomial
volume growth  \cite{Luo2}.

\item Complete biharmonic submanifolds (resp., hypersurfaces) in a negatively curved Riemannian
manifold whose sectional curvature (resp., Ricci curvature) is smaller that $-\epsilon$ for some
$\epsilon>0$  \cite{Luo2}.

\item Complete biharmonic submanifolds (resp., hypersurfaces) $M$ in a Riemannian manifold of non-positive sectional (resp., Ricci) curvature whose mean curvature vector satisfies $\int_M | H^p |dv<\infty$ for some $p>0$  \cite{Luo2}.

\item Complete biharmonic hypersurfaces $M$ in a Riemannian manifold of non-positive Ricci curvature whose mean curvature vector satisfies $\int_M | H|^\alpha dv<\infty$ for some $\epsilon>0$ with $1+\epsilon\leq \alpha<\infty$  \cite{M13}.

\item $\epsilon$-superbiharmonic submanifolds in a complete Riemannian manifolds satisfying the decay condition at infinity
$$\lim_{\rho\to \infty}\frac{1}{\rho^2}\int_{f^{-1}(B_\rho)}| H |^{2}dv=0,$$
where $f$ is the immersion, $B_\rho$ is a geodesic ball of $N$ with radius $\rho$ \cite{Wh}.

\item Proper $\epsilon$-superharmonic submanifolds $M$ with $\epsilon>0$ in a complete Riemannian manifold $N$ whose mean curvature vector satisfying the growth condition
$$\lim_{\rho\to \infty}\frac{1}{\rho^2}\int_{f^{-1}(B_\rho)}| H |^{2+a}dv=0,$$
where $f$ is the immersion, $B_\rho$ is a geodesic ball of $N$ with radius $\rho$, and $a\geq 0$  \cite{Luo2}.
\end{itemize}
\vskip.1in

On the other hand, Ou and  Tang  \cite{OT} proved  that  {\bf generalized biharmonic conjecture is false in general} by constructing foliations of proper biharmonic hyperplanes in some conformally flat 5-manifolds with negative sectional curvature.
Further counterexamples to the generalized biharmonic conjecture were constructed  in \cite{LO}.

\vskip.1in

Finally, the author would like to recall the following two biharmonic conjectures mentioned earlier in  \cite{C32} which are closely related to author's original biharmonic conjecture.
\vskip.1in

{\bf Biharmonic Conjecture for Hypersurfaces}:
\emph{Every biharmonic hypersurface of Euclidean spaces is minimal.}
\vskip.1in

The global version of my original biharmonic conjecture can be found, for instance, in \cite{AM,M13}.
\vskip.1in

 {\bf Global Version of Chen's biharmonic Conjecture}:
\emph{Every complete biharmonic submanifold of a Euclidean space is minimal.}

\end{document}